\newif\ifelsevier
\elsevierfalse

\newif\ifderivations
\derivationsfalse

\ifelsevier
    \documentclass[5p, times]{elsarticle}
    \usepackage{natbib}
    \def\period{}
\else
    \documentclass[12pt]{article}
    \usepackage[margin=1.25in, headheight=42pt]{geometry}
    \usepackage{hyperref}
    \usepackage{fancyhdr}
    \usepackage[style=alphabetic,backend=bibtex,giveninits=true,url=false,doi=false,
                natbib=true,maxbibnames=99]{biblatex}
    \def\period{.}
\fi

\usepackage{float}
\usepackage{listings}
\usepackage{amsmath}
\usepackage{nicefrac}
\usepackage{flushend}
\usepackage{caption}
\usepackage{scalerel,stackengine}

\DeclareCaptionType[fileext=lop]{listing}

\newcommand{\ie}{{\it i.e.}}

\newcommand{\BA}{\begin{array}}
\newcommand{\EA}{\end{array}}

\newcommand{\ones}{\mathbf 1}

\newcommand{\reals}{{\mbox{\bf R}}}

\newcommand{\Expect}{\mathop{\bf E{}}}
\newcommand{\Prob}{\mathop{\bf Prob}}

\newcommand{\argmax}{\mathop{\rm argmax}}

\ifelsevier
\else
    \fancypagestyle{plain}{%
        \fancyhead{}
        \chead{\textsc{For professional clients and qualified investors only \\
               Not for public distribution \\
               This version posted with permission
               }}
    }
\fi

\begin{document}

\ifelsevier
    \begin{frontmatter}
\fi

\title{A Certainty Equivalent Merton Problem}

\ifelsevier
\else
    \author{Nicholas Moehle \and Stephen Boyd}
    \maketitle
\fi

\begin{abstract}
The Merton problem is the well-known stochastic control problem
of choosing consumption over time, as well as an investment mix, to maximize
expected constant relative risk aversion (CRRA)
utility of consumption.  Merton formulated the problem
and provided an analytical solution in 1970; since then a number of
extensions of the original formulation have been solved.
In this note we identify a certainty equivalent problem, \ie,
a deterministic optimal control problem with the same optimal value function
and optimal policy, for the base Merton problem,
as well as a number of extensions.
When time is discretized, the certainty equivalent problem becomes
a second-order cone program (SOCP),
readily formulated and solved
using domain specific languages for convex optimization.
This makes it a good starting point for model predictive control, 
a policy that can handle extensions that are either too cumbersome or
impossible to handle exactly using standard dynamic programming methods.
\end{abstract}

\ifelsevier
        \author{Nicholas Moehle}
        \affiliation{organization={BlackRock AI Labs},
                      %addressline={820 Ramona Ave},
                      city={Palo Alto},
                      %postcode={94301},
                      state={CA},
                      country={USA}}
        \author{Stephen Boyd}
        \affiliation{organization={Electrical Engineering Department, Stanford University},
                      %addressline={David Packard Building, 350 Jane Stanford Way},
                      city={Stanford},
                      %postcode={94305},
                      state={CA},
                      country={USA}}
    \end{frontmatter}
\else
\fi

\section{Introduction}
We revisit Merton's seminal 1970 formulation (and solution) 
of the consumption and investment decisions of an individual investor.
We present a formulation of Merton's problem as a deterministic convex optimal control problem,
and in particular, a second-order cone program (SOCP) when time is discretized.
Even though the Merton problem was first solved more than 50 years ago,
its reformulation as a deterministic
convex optimization problem provides fresh insight into
the solution of the stochastic problem
that may be useful for formulating other multiperiod investment problems
as convex optimization problems.

We also see two practical advantages to the certainty equivalent formulation.
First, for extensions of the Merton problem for which a solution is known,
working out the optimal policy can be complex and error prone.  To handle these 
extensions with the certainty equivalent form, we simply add the appropriate terms 
to the objective or constraints, to obtain the optimal policy.
The problem specification is straightforward and transparent,
especially when expressed 
in a domain specific language (DSL) for convex optimization,
such as \texttt{cvxpy} \citep{cvxpy}.

The second and perhaps more significant advantage is that 
the certainty equivalent problem can be used as a starting point
for further extensions of the Merton problem, for which no closed-form solutions are known.
In this case, the certainty equivalence property is lost,
and solving the deterministic problem no longer solves the corresponding stochastic problem exactly.
We can, however, still use model predictive control (MPC),
a method that involves online convex optimization,
to develop a policy that handles the extension.
MPC policies are simple, easy to implement, fully interpretable,
and have excellent (if not always optimal) practical performance.

\subsection{Previous work}

\paragraph{Merton's problem\period}
Merton's consumption--investment problem dates back to his original
1970 paper \citep{merton1970optimum}.
Many extensions to the basic Merton problem exist,
some of which were covered in Merton's original paper.
(These include deterministic income and general HARA utility.)
Most proposed extensions do not have a closed-form solution,
but some that do include uncertain mortality, life insurance, and annuities,
first adressed by \cite{richard1975optimal}.
Some extensions for the specific case of quadratic utility are handled in \cite{basak2010dynamic}.
We note that many of these extensions individually lead to complicated solutions,
and deriving the optimal policy when several extensions are combined may be very inconvenient
for a practical implementation.
%XXX(more here)
%and Epstein--Zin preferences, first solved in \cite{XXX}.

\paragraph{Certainty equivalence\period}
Rarely, stochastic control problems have a \emph{certainty equivalent} formulation,
\ie, a deterministic optimal control problem with the same optimal policy.
The most famous example is the linear quadratic regulator (LQR) problem,
in which the dynamics are affine, driven by additive noise,
and the stage costs are convex quadratic
\citep{barratt2018stochastic}, \cite[\S 3.1]{bertsekas2017dynamic}, \cite[\S 3]{kwakernaak1972linear}.
In this case, the certainty equivalent problem is obtained by simply ignoring the stochastic noise term.
Many extensions to linear quadratic control also have a certainty equivalent reformulation.
Examples include the linear quadratic Guassian problem, in which the state is imperfectly observed
\citep[\S 5]{kwakernaak1972linear},
and linear exponential quadratic regulator (LEQR) problem, which uses a risk-sensitive cost function
\citep{whittle1990risk}.
Our certainty equivalent formulation is similar to LEQR
in that the uncertain quantity is adversarial 
\citep[\S 10.2]{whittle1990risk}.
(For the Merton problem, the uncertain quantity is the investment returns.)

\paragraph{Model predictive control\period}
In model predictive control, unknown values of future parameters are replaced
with estimates or forecasts over a planning horizon extending from the current time
to some time in the future, resulting in a deterministic optimal control problem.
This problem is solved, with the result intrepretable as a plan of action over the planning horizon.
The MPC policy simply uses the current or first value in the plan of action.
This planning is repeated when updated forecasts are available, using the updated
forecasts and current state.
When applied in the context of stochastic control,
MPC policies are not optimal in general,
but often exhibit excellent practical performance,
and are widely used in several application areas.
MPC is discussed in detail in \citep{borrelli2017predictive, kwon2006receding}.
In \citep{perfbounds}, the authors use a computational bound
to show that MPC is nearly optimal for some stochastic control
problems in finance.

%One justification of model predictive control is its similarity to certainty equivalent LQR.
%We apply the 
As discussed above, ignoring uncertainty is in fact optimal for linear quadratic control,
and MPC leads to an optimal policy when applied to LQR.
In this sense, MPC can be interpreted as applying certainty-equivalence beyond where 
it is theoretically justified
in order to obtain a good heuristic control policy \citep[\S 4.3]{bertsekas2017dynamic}.
For the Merton problem,
we also propose to use a certainty equivalent problem as the basis of an MPC control policy,
even when certain extensions to the base problem ruin exact certainty equivalence.

While MPC has been used in practical applications for decades, recent advances 
make it very attractive, and easy, to develop and deploy.
First, DSLs for convex optimization allow the control policy to be 
expressed in a few lines of very simple and clear code, that express the 
dynamics, objective, and constraints, which makes it easier to develop, debug,
and maintain (for example by adding or updating a constraint).
Code generation systems such as \texttt{cvxgen} 
\citep{cvxgen}
can be used to generate low-level code for that solves the problem specified,
which is suitable for use in high speed embedded applications \citep{wang2009fast}.
In the context of the present paper, this means that the MPC policy we 
propose in \S\ref{s-mpc} can be very conveniently implemented.

\paragraph{Multi-period portfolio optimization\period}
It is instructive to compare our certainty equivalent problem
to popular formulations of multi-period portfolio allocation
(See \cite{cvxportfolio} and references therein).
There are two features present in our certainty equivalent problem
that we do not see in practical multiperiod portfolio construction problems
in the literature:
\begin{enumerate}
\item The risk term (which is quadratic in the dollar-valued asset allocation vector $x_t$),
is normalized by the total wealth $w_t$, which is also is a decision variable.
This risk term is jointly convex in $x_t$ and $w_t$ (and is in fact SOCP representable).
With this normalization,
risk preferences are consistent even as the wealth $w_t$ changes over the investment horizon.
\item The risk term is included as a penalty in the dynamics,
\ie, by taking more risk now, one should expected to have lower wealth in the future.
This contrasts with the tradition of penalizing risk in the objective function.
\end{enumerate}
We believe these to be valuable improvements
to standard multi-period portfolio construction formulations, especially in 
cases when the control or optimization is over a very long time period.

\subsection{Outline}
In \S\ref{s-merton}, we give the base Merton problem and review its solution,
for future reference.
In \S\ref{s-ce-prob}, we give a certainty equivalent problem
and prove equivalence.
In \S\ref{s-extensions}, discuss several extensions to the Merton problem,
and show how each one changes the certainty equivalent formulation.
In \S\ref{s-mpc}, we discuss how to use the certainty equivalent problem 
for model predictive control.

\section{Merton problem}
\label{s-merton}
In this section we discuss the Merton problem and its solution.
To keep the proofs concise, we consider the most basic form of this problem;
extensions are considered in \S\ref{s-extensions}.
Our formulation is in continuous time and relies on stochastic calculus.
However, to maintain both brevity and accessibility,
we are cavalier about the technical details,
with the assumption that a sophisticated reader can fill in the gaps,
or consult other references.

\paragraph{Dynamics\period}
An investor must choose how to invest and consume over
a lifetime of $T$ years.
The investor has wealth $w_t>0$ at time $t$, 
and consumes wealth at rate $c_t>0$, for $t\in [0, T]$,
with the remaining wealth invested in a portfolio with mean rate of return 
$\mu_t$ and volatility $\sigma_t$.
The wealth dynamics are a geometric random walk,
\[
dw_t = (\mu_t w_t - c_t) \, dt + \sigma_t w_t \, dz_t,
\]
where $z_t$ is a Brownian motion.
The initial condition is $w_0 = w_{\rm init} > 0$.

\paragraph{Investment portfolio\period}
The portfolio consists of $n$ assets, with an investment mix given by
the fractional allocation $\theta_t$, with $\ones^T \theta_t=1$ (where $\ones$ is the 
vector with all entries one).  Thus we invest 
$(w_t\theta_t)_i$ dollars in asset $i$,
with a negative value denoting a short position.
The portfolio return rate and volatility are given by
\[
\mu_t = \mu^T \theta_t, \qquad  \sigma_t = (\theta_t^T \Sigma \theta_t)^{1/2},
\]
where $\mu\in \reals^n$ is the mean of the return process, and 
$\Sigma$ is the symmetric positive definite covariance.
(Note that we use the time-varying scalar $\mu_t$ to denote the portfolio
return as a function of time, and the vector $\mu$ to denote the constant expected 
return rates of the $n$ assets.)

The investment allocation decision $\theta_t$ satisfies $\ones^T \theta_t = 1$,
as well as other investment constraints, which we
summarize as $\theta_t \in \Theta$, where $\Theta$ is a convex set.
These could include risk limits, sector exposure limits, or concentration limits.
(See \cite[\S 4.4]{cvxportfolio} for an overview of convex investment constraints.)
For notational convenience, we assume every $\theta_t \in \Theta$ satisfies $\ones^T \theta_t = 1$.

With the portfolio return and volatility we obtain the wealth dynamics
\begin{equation}
\label{e-wealth-dynamics}
dw_t = (\mu^T \theta_t w_t - c_t) \, dt + 
\left(\theta_t^T \Sigma_t \theta_t\right)^{1/2} w_t \, dz_t.
\end{equation}

\paragraph{Utility\period}
The investor has lifetime consumption utility $\int_0^T c_t^\gamma/\gamma \; dt$
and bequest utility $w_T^\gamma / \gamma$.
The risk aversion parameter $\gamma$ satisfies $\gamma < 1$ and $\gamma \neq 0$.
The investor's total expected utility is
\begin{equation}
\label{e-utility}
U = \Expect \bigg(\frac{\beta}{\gamma} w_T^\gamma + \int_0^T \frac{1}{\gamma} c_t^\gamma \, dt \bigg).
\end{equation}
The parameter $\beta>0$ trades off consumption and bequest utility.

%\Expect \bigg(\frac{B}{\gamma} w_T^\gamma 
%            +  \int_t^T \frac{1}{\gamma} c_\tau^\gamma \, d\tau \bigg),
%(The case $\gamma \to 0$ is discussed in \S\ref{s-log-utility}.)
%The coefficient $\alpha_t >0$ is the time weighting of the utility at time $t$.

\paragraph{Stochastic control problem\period}
At each time $t$, the investor chooses the consumption $c_t$ and the investment 
allocation $\theta_t$.
A \emph{policy} maps the time $t$ and the current wealth $w_t$
to the consumption $c_t$ and the allocation $\theta_t$, which we write as
\begin{equation}
\label{e-policy}
(c_t,\theta_t) = \pi_t(w_t),
\end{equation}
where for each $t\in [0,T]$, $\pi_t : \reals_{++}  \to \reals_{++} \times \Theta$.
(Here $\reals_{++}$ denotes the set of positive real numbers.)
The \emph{Merton problem} is to choose a policy $\pi_t$, $t\in [0,T]$,
to maximize $U$.

\subsection{Solution via dynamic programming}
We review here the solution of the Merton problem via dynamic programming,
for completeness and also for future reference.

\paragraph{Value function\period}
The \emph{value function} $V_t: \reals_{++} \to \reals$, for $t\in [0,T]$, is defined as 
\[
V_t(w) = \Expect \bigg(\frac{\beta}{\gamma} w_T^\gamma 
            +  \int_t^T \frac{1}{\gamma} c_\tau^\gamma \, d\tau \bigg),
\]
with $c_\tau$ and $\theta_\tau$ following an optimal policy for $\tau \in [t,T]$,
and initial condition $w_t=w$.
%\[
%V_t(w) = \Expect \bigg(\int_t^T \frac{1}{\gamma} c_t^\gamma \, dt \;\bigg|\; w_t = w \bigg).
%\]
We define $V_T(w) = (B/\gamma) w^\gamma$ for $w>0$.

If the value function is sufficiently smooth, it satisfies the Hamilton-Jacobi-Bellman PDE
\begin{align}\label{e-hjb}
-\dot V_t(w) = \sup_{c, \theta\in\Theta}
\bigg( \frac{1}{\gamma}c^\gamma + V_t'(w) (\mu^T \theta w - c) 
        + \frac12 V_t''(w) (\theta^T \Sigma \theta) w^2 \bigg)
\end{align}
for $w>0$.
Conversely, any function satisfying \eqref{e-hjb} and the terminal condition
$V_T=0$ is the value function.
Here $\dot V_t$ denotes the partial derivative of $V$ with respect to time,
and $V_t'$ and $V_t''$ denote the first and second partial derivatives
with respect to the wealth.

It is well known that the value function for the Merton problem is
\begin{align}\label{e-value-fun-cand}
V_t(w) = a_t \frac{w^\gamma}{\gamma},
\end{align}
where $a_t$ is a function of time.
To obtain $a_t$, we first solve a Markowitz portfolio allocation problem,
%\begin{equation}
%\label{e-markowitz}
%\begin{array}{ll}
%\mbox{maximize} & \displaystyle \mu^T \theta + \frac{\gamma-1}{2} \theta^T \Sigma \theta \\
%\mbox{subject to} & \theta \in \Theta,
%\end{array}
%\end{equation}
\begin{equation}
\label{e-markowitz}
\begin{aligned}
& \mbox{maximize} && \mu^T \theta + \frac{\gamma-1}{2} \theta^T \Sigma \theta \\
& \mbox{subject to} && \theta \in \Theta,
\end{aligned}
\end{equation}
with variable $\theta$.  
(Since $\gamma - 1<0$, the second term is a concave risk adjustment.)
We let $r_{\rm ce}$ denote the optimal value,
and we denote the solution as $\theta_{\rm ce}$.
We then have, for $t\in [0, T]$,
\begin{align}\label{e-b-def}
a_t = 
\bigg(\frac{1-\gamma}{\gamma r_{\rm ce}} 
    \bigg(1 - C \exp\Big(\frac{\gamma r_{\rm ce}}{1-\gamma} (T - t)\Big) \bigg)\bigg)^{1-\gamma},
\end{align}
%where $\alpha = \gamma r_{\rm ce}/(1-\gamma)$.
%
%\begin{align}\label{e-b-def}
%a_t = 
%\bigg(\frac1\alpha \Big(1 - (1 - \alpha B^{1/(1-\gamma)})
%    \exp\big(\alpha(T - t)\big) \Big)\bigg)^{1-\gamma},
%\end{align}
%where $\alpha = \gamma r_{\rm ce}/(1-\gamma)$.
%
%\begin{align}\label{e-b-def}
%a_t = 
%\left(\frac{(1-\gamma) \pig(1 - 
%%\exp \big(\gamma r_{\rm ce}(C + T-t)/(1-\gamma)\big)\pig)}{\gamma r_{\rm ce}}\right)^{1-\gamma},
%C \exp \big(\gamma r_{\rm ce}(T-t)/(1-\gamma)\big)\pig)}{\gamma r_{\rm ce}}\right)^{1-\gamma},
%\end{align}
where $C = 1 - \gamma r_{\rm ce}\beta ^{1/(1-\gamma)} / (1 - \gamma)$.
%\[
%C = \frac{1-\gamma}{\gamma r_{\rm ce}} \log \bigg(1 
%        - \frac{\gamma r_{\rm ce}}{1 - \gamma} \beta ^{1/(1-\gamma)}\bigg).
%\]

\paragraph{Optimal policy\period}
The optimal policy can be expressed in terms of the value function as
\ifelsevier
    \begin{align*}
    \pi_t^\star (w) 
    &= (c_t, \theta_t) \\
    &= \argmax_{c,\theta\in \Theta} \bigg( \frac{1}{\gamma} c^\gamma  
            + V_t'(w) (\mu^T \theta w - c) 
            %\\ & \qquad\qquad\qquad 
            + \frac12 V_t''(w) (\theta^T \Sigma \theta) w^2 \bigg).
    \end{align*}
\else
    \begin{align*}
    \pi_t^\star (w) 
    = (c_t, \theta_t)
    = \argmax_{c,\theta\in \Theta} \bigg( \frac{1}{\gamma} c^\gamma  
            + V_t'(w) (\mu^T \theta w - c) + \frac12 V_t''(w) (\theta^T \Sigma \theta) w^2 \bigg).
    \end{align*}
\fi
With the value function \eqref{e-value-fun-cand}, we obtain the following 
optimal policy.
The consumption has the simple form
\[
c_t = a_t^{1/(\gamma - 1)} w_t,
\]
and the optimal investment mix is constant over time,
\[
\theta_t = \theta_{\rm ce}.
\]
(In extensions of the Merton problem, described below, the
optimal investment mix is not constant over time.)

\paragraph{Proof of optimality\period}
Here we show that the function~\eqref{e-value-fun-cand}
satisfies the Hamilton-Jacobi-Bellman PDE.
To do this, first we substitute $\dot V$, $V_t'$ and $V_t''$ into \eqref{e-hjb}
to obtain
\ifelsevier
    \begin{align*}
    -\dot a_t \frac{w^\gamma}{\gamma} &= \sup_{c, \theta\in\Theta}
    \bigg(
    \frac{1}{\gamma}c^\gamma +
    a_t w^{\gamma-1} (\mu^T \theta w - c)
    \\ & \qquad\qquad
    + \frac12 a_t(\gamma-1) w^{\gamma-2} (\theta^T \Sigma \theta) w^2
    \bigg).
    \end{align*}
\else
    \begin{align*}
    -\overbrace{\dot a_t \frac{w^\gamma}{\gamma}}^{\dot V(w)} &= \sup_{c, \theta\in\Theta}
    \bigg(
    \frac{1}{\gamma}c^\gamma +
    \overbrace{a_t w^{\gamma-1}}^{V'(w)} (\mu^T \theta w - c)
    + \frac12 \overbrace{a_t(\gamma-1) w^{\gamma-2}}^{V''(w)} (\theta^T \Sigma \theta) w^2
    \bigg).
    \end{align*}
\fi
By pulling out $w^{\gamma-1}$ from the last two terms and simplifying,
we obtain
\begin{align}\label{e-hjb-plugged-in}
-\dot a_t \frac{w^\gamma}{\gamma} &=
\sup_{c, \theta\in\Theta} \left(
\frac{1}{\gamma}c^\gamma +
a_t w^{\gamma-1} \left(\Big( \mu^T \theta + \frac{\gamma-1}{2} \theta^T \Sigma \theta \Big) w - c \right) \right).
\end{align}
The maximizing $\theta$ is the solution $\theta_{\rm ce}$ to problem~\eqref{e-markowitz}.
The quantity in the inner parantheses of \eqref{e-hjb-plugged-in}
is the optimal value $r_{\rm ce}$ of this problem,
which can be intrepreted as the certainty equivalent return.
We now have
\begin{align*}
-\dot a_t \frac{w^\gamma}{\gamma} &=
\sup_c \left(\frac{1}{\gamma}c^\gamma + \beta_t w^{\gamma-1} \left(r_{\rm ce} w - c \right) \right).
\end{align*}
%By setting the derivative with respect to $c$ equal to zero and solving, we see that 
The supremum over $c$ is obtained for $c = a_t^{1/(\gamma-1)} w$.
%\begin{align*}
%c^{\gamma-1} = a_t w^{\gamma-1}.
%\end{align*}
Substituting in this value and simplifying, we obtain
\ifderivations
    \begin{align*}
    -\dot a_t \frac{w^\gamma}{\gamma} &=
    \frac{1}{\gamma}a_t^{\gamma/(\gamma-1)}w^\gamma 
        + a_t w^{\gamma-1} \left(r_{\rm ce} w - a_t^{1/(\gamma-1)}w \right)
    \end{align*}
    or,
    \begin{align*}
    -\dot a_t \frac{1}{\gamma} 
    &=
    \frac{1}{\gamma}a_t^{\gamma/(\gamma-1)} + a_t r_{\rm ce} - a_t^{1/(\gamma-1) + 1}\\
    &=
    \frac{1}{\gamma}a_t^{\gamma/(\gamma-1)} + a_t r_{\rm ce} - a_t^{\gamma/(\gamma-1)}\\
    &=
    \left(\frac{1}{\gamma} - 1\right) a_t^{\gamma/(\gamma-1)} + a_t r_{\rm ce} \\
    &=
    \frac{1 - \gamma}{\gamma}a_t^{\gamma/(\gamma-1)} + a_t r_{\rm ce}
    \end{align*}
    or,
\fi
    \begin{align*}
    -\dot a_t = (1 - \gamma)a_t^{\gamma/(\gamma-1)} + \gamma a_t r_{\rm ce}.
    \end{align*}
It can be verified that the definition of $a_t$
in \eqref{e-b-def}
is indeed a solution to this differential equation with terminal 
condition $a_T = \beta$.

\ifderivations
\paragraph{\it Derivation\period}
Define $\alpha_t$ such that
\[
a_t = -\alpha_t^{1-\gamma}.
\]
and $\eta_t$ such that
\begin{align*}
\eta_t &= \frac{\exp \big(\gamma r_{\rm ce}(T-t)/\delta\big)}{\delta} \\
&= 1 + \frac{\gamma r_{\rm ce}}{\delta} a_t^{\frac{1}{1-\gamma}}.
\end{align*}
This means
From the definition of $a_t$, we have
\[
\dot a_t = (1-\gamma)\alpha^{-\gamma}
\underbrace{\frac{\exp \big(\gamma r_{\rm ce}(T-t)/\delta\big)}{\delta}}_{\eta_t}.
\]
With this definition of $\eta_t$, we have
\[
\dot a_t = (1-\gamma)a_t^{\frac{\gamma}{\gamma-1}} + \gamma r_{\rm ce} a_t.
\]
\fi

\section{Certainty equivalent problem}
\label{s-ce-prob}
In this section we present a deterministic convex optimal control problem
that is equivalent to the Merton problem
in the sense that it has the same value function and same optimal policy.

This certainty equivalent problem is
%\begin{equation}
%\label{e-ce-prob}
%\begin{array}{ll} \mbox{maximize} & \displaystyle \int_0^T \frac{1}{\gamma} c_t^\gamma \, dt \\
%\mbox{subject to} & \displaystyle \dot w_t \leq \mu^T x_t - c_t 
%                    + \frac{(\gamma - 1)}{2} \frac{x_t^T \Sigma x_t}{w_t},
%                    \quad t\in [0, T]\\
%                  & w_t = \mathbf 1^T x_t, \quad t\in [0, T]\\
%                  & x_t /w_t \in \Theta, \quad t\in [0,T]\\
%                  & w_0 = w_{\rm init}.
%\end{array}
%\end{equation}
\begin{equation}
\label{e-ce-prob}
\begin{aligned}
& \mbox{maximize} &&
    \frac{\beta}{\gamma} w_T^\gamma + \int_0^T \frac{1}{\gamma} c_t^\gamma \, dt \\
& \mbox{subject to}
  && \dot w_t \leq \mu^T x_t - c_t + \frac{(\gamma - 1)}{2} \frac{x_t^T \Sigma x_t}{w_t}, \quad t\in [0, T]\\
%& && w_t = \mathbf 1^T x_t, \quad t\in [0, T]\\
& && x_t /w_t \in \Theta, \quad t\in [0,T]\\
& && w_0 = w_{\rm init}.
\end{aligned}
\end{equation}
The variables are the consumption $c_t:[0,T]\to\reals_{++}$,
wealth $w_t:[0,T]\to\reals_{++}$,
and $x_t:[0, T]\to\reals^n$,
which is the dollar-valued allocation of wealth to each asset.
(In the notation of \S\ref{s-merton}, we have $x_t = w_t \theta_t$, and $\theta_t = x_t/w_t$.)
%We use the convention that 
%for $w_t = 0$, we have $x_t/w_t \in \Theta$ if and only if $x_t = 0$.
%Furthermore, for $x_t = 0$ and $w_t = 0$, we define $x_t^T \Sigma x_t / w_t=0$.
Note that the constraint $x_t / w_t \in \Theta$ implies 
%$\ones^T (x_t / w_t) = 1$, or equivalently, 
$\ones^T x_t = w_t$,
\ie, the total wealth is the sum of the dollar-valued asset allocations.
%However, we leave the explicit constraint $1^T x_t = w_t$ in \eqref{e-ce-prob} for clarity.

The objective is the lifetime utility, but without expectation since 
this problem is deterministic.
The first constraint resembles the dynamics
of the stochastic process~\eqref{e-wealth-dynamics},
and we call this the \emph{dynamics constraint}.
We will see that for any solution to \eqref{e-ce-prob},
this inequality constraint holds with equality,
in which case the dynamics constraint becomes a (deterministic) ODE.

\paragraph{Interpretation\period}
The problem can be interpreted in the following way.
We plan for a single outcome of the stochastic process~\eqref{e-wealth-dynamics}.
In particular, the dynamics constraint restricts
the growth rate of the wealth to be no greater than the $\mu^T x_t - c_t$
(the mean growth rate in the stochastic process~\eqref{e-wealth-dynamics}),
but reduced by the additional term $(1/2)(\gamma-1) x_t^T \Sigma x_t / w_t$.
Because $\gamma < 1$, this term is negative.
With the change of  variables $\theta_t = x_t / w_t$,
we have
\[
\frac{x_t^T \Sigma x_t}{w_t} = w_t \theta_t^T \Sigma \theta_t,
\]
\ie, this adjustment term is proportional to the variance of the portfolio growth rate
with investment allocation $\theta_t = x_t / w_t$.
In other words, we are pessimistically planning for bad investment returns,
with the degree of pessimism depending on the risk aversion parameter
$\gamma$ and the risk of our portfolio.

In fact, in problem~\eqref{e-ce-prob}, we plan for the returns
\[
r_t = \mu + \frac{\gamma - 1}{2w_t} \Sigma x_t =
\mu + \frac{\gamma - 1}{2} \Sigma \theta_t.
\]
The coefficients in front of $\Sigma x_t$ and $\Sigma\theta_t$ are negative, and 
the entries of $\Sigma x_t$ and $\Sigma \theta_t$ are typically positive.
The vector $\Sigma \theta_t$ can be interpreted as the risk allocation
to the individual assets in the portfolio, since
\[
\theta_t^T \Sigma \theta_t = \sum_{i=1}^n (\theta_t)_i \left(\Sigma \theta_t\right)_i.
\]
In other words, the planned asset returns are the mean returns,
reduced in proportion to
the marginal contribution of each asset to the portfolio variance.
This is related to the concept of risk parity \citep{bai2016least}.

\paragraph{Convexity\period}
Convexity of \eqref{e-ce-prob} follows from the fact that the risk penalty term
$x_t^T \Sigma x_t / w_t$ is a quadratic-over-linear function,
with is jointly convex in $x_t$ and $w_t$ \cite[\S 3.1.5]{cvxbook}.
Also, the set 
\[
\{(x_t, w_t) \in\reals^{n} \times \reals_{++} \mid x_t/w_t \in \Theta \}
\]
is the perspective of $\Theta$,
which is convex when $\Theta$ is \citep[\S 2.3.3]{cvxbook}.
In fact, in most practical portfolio construction problems,
$\Theta$ can described by a collection of linear and quadratic constraints
\citep[\S 4.4]{cvxportfolio}.
In this case, when problem~\eqref{e-ce-prob} is discretized, it becomes an SOCP,
which we describe in \S\ref{s-discretized-prob}.

\paragraph{Equivalence to Merton problem\period}
The Merton problem and problem~\eqref{e-ce-prob} are equivalent
in the sense that they have the same value function
and optimal policy.

To see this, we first consider a modified version of \eqref{e-ce-prob}
in which we convert the dynamics to an equality constraint using a slack variable $u_t \ge 0$:
\[
\dot w_t = \mu^T x_t - c_t + \frac{(\gamma - 1)}{2} \frac{x_t^T \Sigma x_t}{w_t} + u_t.
\]
The new control input $u_t$ can be interpreted as the rate at which we discard wealth.
(We will see that at optimality $u_t = 0$.)
For this modified problem, the Hamilton-Jacobi-Bellman equation is
\begin{align*}
-\dot V(w) &= \sup_{c, x\in w\Theta, u \ge 0}
\frac{1}{\gamma}c^\gamma +
V_t'(w) \left(\Big( \mu^T x + \frac{\gamma-1}{2w} x^T \Sigma x \Big) w - c - u \right).
\end{align*}
First note that with our value function candidate~\eqref{e-value-fun-cand},
we have $V'(w) > 0$, and therefore $u = 0$, as expected.
Now, by using the change of variables $x = \theta w$
and plugging in our value function candidate,
this equation becomes \eqref{e-hjb-plugged-in}.
From this point on,
the proof that this candidate value function satisfies the 
Hamilton-Jacobi-Bellman equation
proceeds exactly as for the (stochastic) Merton problem.

\section{Exact extensions}
\label{s-extensions}
Here we consider several extensions to the Merton problem,
all of which are known in the literature and have closed-form solutions.
For each one, we describe how to modify problem \eqref{e-ce-prob}
to maintain the certainty-equivalence property.

%\paragraph{Discounted utility\period}
%If we replace the utility with
%\[
%U = \int_0^T \frac{\alpha_t}{\gamma} c^\gamma,
%\]
%where $\alpha_t >0$ is the discount of the utility at time $t$,
%then the objective of the certainty equivalent problem becomes
%\[
%\int_0^T \frac{\alpha_t}{\gamma} c^\gamma.
%\]

\paragraph{Time-varying parameters\period}
The Merton problem can be solved when $\mu$, $\Sigma$, and $\Theta$ change over time.
To handle this in the certainty equivalent problem,
we simply replace these parameters by $\mu_t$, $\Sigma_t$, and $\Theta_t$.
(Here $\mu_t$ denotes the time-varying vector of asset expected returns,
a notation clash with our previous use of $\mu_t$ as the scalar portfolio 
expected return.)
Similarly, if we discount the consumption utility of the Merton problem:
\begin{equation*}
U = \Expect \bigg(\frac{\beta}{\gamma} w_T^\gamma + 
        \int_0^T \frac{\alpha_t }{\gamma} c_t^\gamma \, dt \bigg).
\end{equation*}
where $\alpha_t >0$ is the discount of the consumption utility at time $t$,
then the objective of the certainty equivalent problem will match $U$
(but without the expectation).

\paragraph{Uncertain mortality and bequest\period}
Here the terminal time $t_f \in [0, T]$ is random with probability density $p_t$
and \emph{survival function}
\[
s_t = \Prob(t_f > t) = \int_t^T p_t \, dt.
\]
In this case, the investor's utility is
\[
U = \Expect \bigg( \frac{B}{\gamma} w_{t_f}^\gamma 
        + \int_0^{t_f} \frac{1}{\gamma} c_t^\gamma \, dt \bigg).
\]
Here the expectation is taken over $t_f$ as well as 
the paths of the stochastic process~\eqref{e-wealth-dynamics}.

With this modification, the objective of the certainty equivalent problem changes to
\[
\int_0^T \left( \frac{p_t \beta}{\gamma} w_t^\gamma + \frac{s_t}{\gamma} c_t^\gamma \right) \, dt.
\]
We weight the consumption utility by the probability the investor is still alive,
\ie, we treat the survival function as a discount factor.
We also get utility for the bequest continuously over the interval $[0, T]$,
weighted by the density function $p_t$.

\paragraph{Annuities and life insurance\period}
This extension is due to \cite{richard1975optimal}.
Continuing with the previous extension,
we allow the investor to purchase life insurance.
The premium is $l_t$, which the investor can choose,
and the payout of the plan is $\lambda_t l_t$,
where $\lambda_t \ge 0$ is the payout-to-premium ratio at time $t$.
When $l_t<0$, we interpret this as an annuity.
In particular, at time $t$, the investor has $-l_t$ in the annuity account,
which is lost on death,
in return for an additional return of $-\lambda_t l_t$.
The actuarially fair value of $\lambda_t$ is $p_t / s_t$,
which is called the \emph{force of mortality}.
(If $\lambda_t > p_t/s_t$, then life insurance is favorable and annuities are unfavorable;
if $\lambda_t < p_t/s_t$, the reverse is true.) 

With this modification,
the objective of the certainty-equivalant problem changes to
\[
U = \int_0^T \left( \frac{p_t \beta}{\gamma} (w_t + \lambda_t l_t)^\gamma + \frac{s_t}{\gamma} c_t^\gamma \right) \, dt,
\]
\ie, we add the insurance payout to the wealth in the bequest utility.
The dynamics change to
\[
\dot w_t \leq \mu^T x_t - c_t - l_t + \frac{(\gamma - 1)}{2} \frac{x_t^T \Sigma x_t}{w_t}.
\]
Here we subtract the insurance premium from the growth rate of the wealth.

\paragraph{Income\period}
We can add a deterministic income stream, with income rate $y_t$ at time $t$.
The stochastic dynamics are modified be the addition of $y_t$ to the drift term of the wealth process, \ie,
\[
\mu_t = \mu^T \theta_t w_t + y_t - c_t.
\]
In this case, we also assume one of the assets is risk free
with return $\mu_{\rm rf}$ and volatility $0$,
and that 
\begin{equation}
\label{e-simple-theta}
\Theta = \{\theta \mid \ones^T \theta = 1\}.
\end{equation}
These assumptions
allow the investor to counteract the income stream by shorting the risk-free asset
and investing the proceeds in a preferred portfolio of other assets.
The fair value of the income stream is its net present value over $[t, T]$
at the risk-free rate:
\[
v_t = \int_t^T e^{-\mu_{\rm rf} y_t} \, dt,
\]
which can be interpreted as the remaining \emph{human capital} of the investor.
%The value function of the Merton problem is modified to be
%\[
%V_t(w) = a_t \frac{(w + v_t)^\gamma}{\gamma},
%\]
%\ie, the wealth variable is augmented by the total human capital at time $t$.
%The derivation of the certainty equivalent problem proceeds similarly to the standard case;

For this extension, the dynamics in \eqref{e-ce-prob} are replaced by
\[
\dot w_t \leq \mu_t^T x_t + y_t - c_t + \frac{(\gamma - 1)}{2} \frac{x_t^T \Sigma x_t}{w_t + v_t}.
\]
Note the addition of the income term $y_t$
and the normalization of risk by the total wealth plus the remaining human capital.
In this case, the wealth $w_t$ need not be positive
but instead satisfies $w_t + v_t > 0$.
Because of this, we also replace the constraint $x_t / w_t \in \Theta$ 
(which is not defined for $w_t = 0$) with $\ones^T x_t = w_t$.

\paragraph{Epstein--Zin preferences\period}
One interesting feature of the certainty equivalent problem~\eqref{e-ce-prob} is that the risk aversion
parameter $\gamma$ appears separately in the objective and dynamics constraint.
It is reasonable to ask whether, by modifying the consumption utility to be
\[
%\int_0^T \frac{1}{\rho} c_t^\rho \, dt,
\frac{\beta}{\rho} w_T^\rho + \int_0^T \frac{1}{\rho} c_t^\rho \, dt
\]
for some $\rho \neq \gamma$ with $\rho < 1$ and $\rho \neq 0$,
but keeping $\gamma$ in the dynamics constraint,
problem~\eqref{e-ce-prob} is equivalent to some variant of the Merton problem.
This is indeed the case,
but with the expected utility $U$ replaced by Epstein--Zin preferences,
where $1/\rho$ is the elasticity of intertemporal substitution
and $\gamma$ is the risk aversion.
%(For  \cite{duffie1992stochastic}),
%An exact definition of Epstein--Zin preferences in continuous time is complicated;
For details, see \cite{duffie1992stochastic}.

\section{Inexact extensions}
\label{s-inexact-extensions}
Here we discuss several extensions of problem \eqref{e-ce-prob}
that (to our knowledge) do not exactly solve a version of the Merton problem.
Some of these build on the exact extensions of \S\ref{s-extensions}.

\paragraph{Modified utility\period} 
We can change the objective of \eqref{e-ce-prob} to use any
increasing, concave utility function for either consumption or bequest.
These utility functions need not be additive over time:
For example, we can maximize the minimum consumption over the interval $[0, T]$,

As a special case, we can add a minimum consumption constraint
\[
c_t \ge c_t^{\rm min},
\]
where $c_t^{\rm min}$ is the minimum allowable consumption amount
as a function of age.
Similarly, we can enforce a minimum bequest over some time window
(say, to care for underage dependents until they come of age).

\paragraph{Spending limit\period}
We can limit consumption as a fraction of income with the constraint
\[
c_t \le \eta y_t
\]
for some parameter $\eta > 0$.
For example, when $\eta = 0.7$,
this constraint means that we can't consume more than 70\% of our income,
\ie, we must have a savings rate of 30\%.

This constraint can be adjusted to account for investment income.
To see this, take $d\in\reals^n$ to be the vector of dividend yields for each asset,
which is constant and known in advance.
The modified constraint becomes
\[
c_t \le \eta y_t + d^T x_t.
\]
When this constraint is tight, \ie, when we desire to consume more than $\eta$ times our income,
there is added incentive to invest in assets with high dividend yield.

\paragraph{Minimum cash balance\period}
We can include a constraint that the amount invested in cash
be above a certain level,
\ie,
\[
(x_t)_i \ge (x_t^{\rm min})_i,
\]
where $i$ is the index of the cash asset.
This is similar to an emergency fund constraint
that we must keep six months worth of consumption in cash,
which is expressed as
\[
(x_t)_i \ge 0.5 c_t.
\]

\section{Application to model predictive control}
\label{s-mpc}
Model predictive control is a technique for stochastic control problems
that leverages a deterministic approximation of the stochastic problem.
To evaluate an MPC policy, we first solve this determistic problem
to obtain a planned trajectory for the state and control input over the planning horizon.
We then implement only the first control input in this plan,
and rest of the planned trajectory is discarded.
To obtain future control inputs,
the policy is evaluated again,
which requires solving a new deterministic problem.
%If the MPC policy is evaluated again, the process is repeated.

In the context of the Merton problem,
the certainty equivalent problem is used as a basis for a simple model predictive control policy,
which we denote $\pi_t^{\rm mpc}$.
We first define this policy when $t = 0$, with initial wealth $w_0$.
We start by solving the deterministic control problem \eqref{e-ce-prob}
to obtain the optimal trajectories $c_t$ and $\theta_t$.
The MPC policy then takes $\pi_0^{\rm mpc}(w_0) = (c_0,\theta_0)$.
%to the control inputs $c_0$ and $\theta_0$ are optimal 
To define the MPC policy for $t \in (0,T)$,
we first form a new instance of problem~\eqref{e-ce-prob},
which is defined over the interval $[t, T]$ and has initial wealth $w_t$.
Once again we solve the deterministic optimal control problem \eqref{e-ce-prob},
to obtain optimal $c_\tau$ and $\theta_\tau$ over the interval
$\tau \in [t,T]$.
We then take $\pi_t^{\rm mpc}(w_t) = (c_t,\theta_t)$.
Evaluating the MPC policy therefore always requires solving a deterministic 
optimal control problem of the form \eqref{e-ce-prob}.

MPC is a convenient way to implement the optimal policy
for the basic problem or any of the extensions of \S\ref{s-extensions}.
In those cases, the MPC policy is optimal.
When MPC is applied with constraints and an objective
that do not correspond to any version of Merton problem,
the MPC policy is a sophisticated heuristic, and very useful in practice.
%In this sense, the certainty equivalent result is only a launching point
%to develop interesting convex-optimization-based planning algorithms.

To use MPC in practice requires discretizing problem~\eqref{e-ce-prob},
which we discuss in the next section.

\section{Discretized problem}
\label{s-discretized-prob}
Here we show how to discretize problem~\eqref{e-ce-prob}.
We do this for the basic problem only,
but note that the extensions can be handled similarly.

We let $x_k$ denote the value of $x_t$ in \eqref{e-ce-prob} at time
$t = hk$, $k=0,\ldots, K$, where $h=T/K$ is the discretization interval.
(We use the same notation, but index $x$ with the subscript $k$
to denote the discretized variable, and index with $t$ to denote the 
continuous variable.)
We similarly define the discretized variables $c_k$ and $w_k$.
Replacing the time derivative $\dot w_t$ with the forward Euler approximation
$(w_{k+1}-w_k)/h$, and replacing the integral in the objective with a 
Riemann sum approximation, we obtain the discretized problem
%\begin{equation}
%\label{e-ce-prob-discrete}
%\begin{array}{ll} \mbox{maximize} & (1/\gamma) \sum_{k=0}^{K} c_k^\gamma \\
%\mbox{subject to} & \displaystyle \frac{w_{k+1} - w_k}{h} \le \mu^T x_k - c_k 
%                    + \frac{(\gamma - 1)}{2} \frac{x_k^T \Sigma x_k}{w_k},
%                    %+ (1/2)(\gamma - 1) x_k^T \Sigma x_k/w_k,
%                    \quad k=0,\dots, K-1 \\
%                  & w_k = \mathbf 1^T x_k, \quad k=0,\dots, K \\
%                  & x_k /w_k \in \Theta, \quad k=0,\dots, K \\
%                  & w_0 = w_{\rm init}
%                  %& w_K \ge 0 \\
%\end{array}
%\end{equation}
\ifelsevier
    \begin{equation}
    \label{e-ce-prob-discrete}
    \begin{aligned}
    & \mbox{maximize}   && \frac{\beta}{\gamma} w_T^\gamma + \sum_{k=0}^{K-1} \frac h\gamma c_k^\gamma \\
    & \mbox{subject to} && \displaystyle \frac{w_{k+1} - w_k}{h} \le \mu^T x_k - c_k 
                            + \frac{(\gamma - 1)}{2} \frac{x_k^T \Sigma x_k}{w_k} \\
    %&                   && w_k = \mathbf 1^T x_k \\
    &                   && x_k /w_k \in \Theta \\
    &                   && w_0 = w_{\rm init}.
    \end{aligned}
    \end{equation}
\else
    \begin{equation}
    \label{e-ce-prob-discrete}
    \begin{aligned}
    & \mbox{maximize}   && \frac{\beta}{\gamma} w_T^\gamma + \sum_{k=0}^{K-1} \frac h\gamma c_k^\gamma \\
    & \mbox{subject to} && \displaystyle \frac{w_{k+1} - w_k}{h} \le \mu^T x_k - c_k 
                            + \frac{(\gamma - 1)}{2} \frac{x_k^T \Sigma x_k}{w_k},
                            \quad k=0,\dots, K-1 \\
    %&                   && w_k = \mathbf 1^T x_k, \quad k=0,\dots, K \\
    &                   && x_k /w_k \in \Theta, \quad k=0,\dots, K \\
    &                   && w_0 = w_{\rm init}.
    \end{aligned}
    \end{equation}
\fi
The variables are $x_k\in\reals^n$ and $w_k\in\reals_{++}$ for $k=0,\dots,K$
and $c_k\in\reals_{++}$ for $k=0,\dots,K-1$.
\ifelsevier
    The first constraint holds for $k = 0, \dots, K-1$,
    and the second constraint holds for $k = 0, \dots, K$.
\fi
%Here $0, \dots, K$ are the discrete time periods representing the time interval $[0, T]$,
%and $h$ is the time step parameter chosen such that $T = Kh$.
All of the extensions (exact and inexact) discussed above
can be discretized as well, but we do not give the details here.

The discretized certainty equivalent
problem~\eqref{e-ce-prob-discrete} is a (finite-dimensional)
convex optimization problem, and can therefore be 
easily expressed in a domain-specific language for convex optimization,
such as \texttt{cvxpy}.
As an example, 
we give a \texttt{cvxpy} implementation of \eqref{e-ce-prob-discrete}
in listing~\ref{l-cvxpy-prob}
when $\Theta$ is given by \eqref{e-simple-theta}.

\lstset{frame=lines}
\lstset{basicstyle=\small\ttfamily,breaklines=true}
\begin{listing}
\begin{lstlisting}
w = Variable(K+1)
x = Variable(n, K+1)
c = Variable(K)
Sigma_half = numpy.linalg.cholesky(Sigma)

U = beta/gamma * power(w[K], gamma) + h/gamma * sum(power(c, gamma))
constr = [w == sum(x, axis=0), w[0] == w_init]
for k in range(K):
    constr += [diff(w[k+1] - w[k])/h <= mu @ x[:, k] - c[k] + (gamma - 1)/2 * quad_over_lin(Sigma_half @ x[:, k], w[k])]
]
problem = Problem(Maximize(U), constr)
problem.solve()
\end{lstlisting}
\caption{
    An implementation of the discretized certainty
    equivalent problem~\eqref{e-ce-prob-discrete} using \texttt{cvxpy}.
}
\label{l-cvxpy-prob}
\end{listing}

For most practical portfolio construction problems,
$\Theta$ is SOCP representable,
which means that problem~\eqref{e-ce-prob-discrete} is an SOCP \citep{socp}.
To see this, note that the power utility $c_k^\gamma$ 
and the quadratic-over-linear functions are SOCP representable; see 
\cite[\S 2.2.f]{alizadeh2003second} and \cite[\S 2.4]{socp}, respectively.
The perspective of $\Theta$
can be represented using the same cones used to represent $\Theta$
\citep[\S 2]{moehle2015perspective}.

To give some idea of the speed at which current solvers
can solve the discretized problem \eqref{e-ce-prob-discrete}
(and its extensions), consider a problem with $n=500$ assets, $K=50$ periods,
and covariance matrix $\Sigma$ given as a typical factor model, with 25 factors.  
This problem has more than 100000 optimization variables.
With just a small modification 
of the code given in listing \ref{l-cvxpy-prob} to exploit the
low rank plus diagonal structure of the covariance matrix, 
the open-source solver ECOS \cite{ecos} solve the problem in around 
two seconds, on a single thread.

\ifelsevier
    \bibliographystyle{plainnat}
    \bibliography{merton}

@book{cvxbook,
  title={Convex optimization},
  author={S. Boyd and L. Vandenberghe},
  year={2004},
  publisher={Cambridge University Press}
}

@article{bai2016least,
  title={Least-squares approach to risk parity in portfolio selection},
  author={Bai, Xi and Scheinberg, Katya and Tutuncu, Reha},
  journal={Quantitative Finance},
  volume={16},
  number={3},
  pages={357--376},
  year={2016},
  publisher={Taylor \& Francis}
}

@article{alizadeh2003second,
  title={Second-order cone programming},
  author={Alizadeh, Farid and Goldfarb, Donald},
  journal={Mathematical programming},
  volume={95},
  number={1},
  pages={3--51},
  year={2003},
}

@article{moehle2015perspective,
  title={A perspective-based convex relaxation for switched-affine optimal control},
  author={Moehle, Nicholas and Boyd, Stephen},
  journal={Systems and Control Letters},
  volume={86},
  pages={34--40},
  year={2015},
  publisher={Elsevier}
}

@inproceedings{ecos,
  author={Domahidi, Alexander and Chu, Eric and Boyd, Stephen},
  booktitle={European Control Conference},
  title={{ECOS}: {A}n {SOCP} solver for embedded systems},
  year={2013},
  pages={3071-3076}
}

@article{wang2009fast,
  title={Fast model predictive control using online optimization},
  author={Wang, Yang and Boyd, Stephen},
  journal={IEEE Transactions on control systems technology},
  volume={18},
  number={2},
  pages={267--278},
  year={2009},
  publisher={IEEE}
}

@article{cvxgen,
  title={{CVXGEN}: A code generator for embedded convex optimization},
  author={Mattingley, Jacob and Boyd, Stephen},
  journal={Optimization and Engineering},
  volume={13},
  number={1},
  pages={1--27},
  year={2012},
  publisher={Springer}
}

@unpublished{barratt2018stochastic,
  title={Stochastic control with affine dynamics and extended quadratic costs},
  author={Barratt, Shane and Boyd, Stephen},
  note={ArXiv preprint},
  year={2018}
}

@article{socp,
  title={Applications of second-order cone programming},
  author={Lobo, Miguel Sousa and Vandenberghe, Lieven and Boyd, Stephen and Lebret, Herv{\'e}},
  journal={Linear algebra and its applications},
  volume={284},
  number={1-3},
  pages={193--228},
  year={1998},
  publisher={Elsevier}
}

@article{cvxpy,
  author  = {S. Diamond and S. Boyd},
  title   = {{CVXPY}: {A} {P}ython-embedded modeling language for convex optimization},
  journal = {Journal of Machine Learning Research},
  year    = {2016},
  volume  = {17},
  number  = {83},
  pages   = {1--5},
}

@book{borrelli2017predictive,
  title={Predictive control for linear and hybrid systems},
  author={Borrelli, Francesco and Bemporad, Alberto and Morari, Manfred},
  year={2017},
  publisher={Cambridge University Press}
}

@book{kwon2006receding,
  title={Receding horizon control: Model predictive control for state models},
  author={Kwon, Wook Hyun and Han, Soo Hee},
  year={2006},
  publisher={Springer}
}

@book{whittle1990risk,
  title={Risk-sensitive Optimal Control},
  author={Whittle, Peter},
  year={1990},
  publisher={John Wiley \& Sons}
}

@book{bertsekas2017dynamic,
  title={Dynamic programming and optimal control},
  author={Bertsekas, Dimitri P},
  volume={4},
  year={2017},
  publisher={Athena scientific}
}

@incollection{merton1970optimum,
  title={Optimum consumption and portfolio rules in a continuous-time model},
  author={Merton, Robert C},
  booktitle={Stochastic Optimization Models in Finance},
  pages={621--661},
  year={1970},
  publisher={Elsevier}
}

@book{kwakernaak1972linear,
  title={Linear optimal control systems},
  author={Kwakernaak, Huibert and Sivan, Raphael},
  volume={1},
  year={1972},
  publisher={John Wiley \& Sons}
}

@article{cvxportfolio,
  title={Multi-period trading via convex optimization},
  author={Boyd, Stephen and Busseti, Enzo and Diamond, Steve and Kahn, Ronald N and Koh, Kwangmoo and Nystrup, Peter and Speth, Jan},
  journal={Foundations and Trends in Optimization},
  volume={3},
  number={1},
  pages={1--76},
  year={2017},
  publisher={Now Publishers, Inc.}
}

@article{perfbounds,
  title={Performance bounds and suboptimal policies for multi-period investment},
  author={Boyd, Stephen and Mueller, Mark and O'Donoghue, Brendon and Wang, Yang},
  journal={Foundations and Trends in Optimization},
  volume={1},
  number={1},
  pages={1--69},
  year={2014},
  publisher={Now Publishers, Inc.}
}

@article{duffie1992stochastic,
  title={Stochastic differential utility},
  author={Duffie, Darrell and Epstein, Larry G},
  journal={Econometrica: Journal of the Econometric Society},
  pages={353--394},
  year={1992},
  publisher={JSTOR}
}

@article{richard1975optimal,
  title={Optimal consumption, portfolio and life insurance rules for an uncertain lived individual in a continuous time model},
  author={Richard, Scott F},
  journal={Journal of Financial Economics},
  volume={2},
  number={2},
  pages={187--203},
  year={1975},
  publisher={Elsevier}
}

@article{basak2010dynamic,
  title={Dynamic mean-variance asset allocation},
  author={Basak, Suleyman and Chabakauri, Georgy},
  journal={The Review of Financial Studies},
  volume={23},
  number={8},
  pages={2970--3016},
  year={2010},
  publisher={Oxford University Press}
}
\else
    \clearpage
    \printbibliography
\fi

\end{document}